\documentclass{article}

\begin{document}

\title{Some recent developments in K\"ahler geometry and exceptional holonomy}
\author{Simon Donaldson\\ Simons Centre for Geometry and Physics, Stony Brook\\ \\Imperial College,
  London}
  

\maketitle


\newcommand{\bR}{{\bf R}}
\newcommand{\bC}{{\bf C}}
\newcommand{\bZ}{{\bf Z}}
\newcommand{\bP}{{\bf P}}
\newcommand{\bH}{{\bf H}}
\newcommand{\dbd}{\partial\overline{\partial}}
\section{Introduction}
This article is a broad-brush survey of two areas in differential geometry. While these two areas are  not usually put side-by-side in this way, there are several reasons for   discussing them together. First, they both fit into a very general pattern, where one asks about the existence  of various differential-geometric structures on a manifold. In one case we consider a complex K\"ahler manifold and seek a distinguished metric, for example a K\"ahler-Einstein metric. In the other we seek a metric of exceptional holonomy on a manifold of dimension $7$ or $8$. Second, as we shall see in more detail below, there are numerous points of contact between these areas at a  technical level. Third, there is a pleasant contrast  between the state of development in the  fields. These questions in K\"ahler geometry have been studied for more than half a century: there is a huge literature with many deep and far-ranging results. By contrast, the theory of manifolds of exceptional holonomy is a wide-open field: very little is known in the way of general results and the developments so far have focused on examples. In many cases these examples depend on advances in K\"ahler geometry.

\section{K\"ahler geometry}
\subsection{Review of basics}

Let $X$ be a compact complex manifold of complex dimension $n$. We recall that a Hermitian metric on $X$ corresponds to a positive 2-form $\omega_{0}$ of type $(1,1)$ That is, in local complex co-ordinates $z_{a}$
$$\omega_{0}= i \sum_{a, b} g_{ab} dz_{a}d\overline{z}_{b}, $$
where at each point $\left(g_{ab}\right) $ is a positive-definite Hermitian matrix. The Hermitian metric is {\it K\"ahler} if $\omega_{0}$ is a closed form, so then $X$ is simultaneously a complex, Riemannian and symplectic manifold.   The form $\omega_{0}$ defines a cohomology class $[\omega_{0}]\in H^{2}(X;\bR)$ and any other K\"ahler metric in this class can be represented by a potential function $\omega_{\phi}=\omega_{0}+ i \dbd \phi$. In local co-ordinates this changes $g_{ab}$ to
$$   g_{ab}+ \frac{\partial^{2}\phi}{\partial z_{a}\partial \overline{z}_{b}} . $$

Recall that in Riemannian geometry the {\it Ricci tensor} is a symmetric 2-tensor given by a contraction of the full Riemann curvature. Geometrically, it represents the infinitesimal behaviour of the volume form of the manifold, in geodesic co-ordinates. In the K\"ahler case the Ricci curvature also corresponds to a $(1,1)$ form and has another interpretation as  the curvature of the induced connection on the {\it anti-canonical line bundle} $K^{-1}_{X}=\Lambda^{n}TX$. In local co-ordinates this Ricci form is given by
$   \rho= i\dbd L, $
where $L$ is the logarithm of the volume element
$$   L =\log \det \left(  g_{ab}+ \frac{\partial^{2}\phi}{\partial z_{a}\partial \overline{z}_{b}}\right). $$
The trace of the Ricci curvature is the {\it scalar curvature} $S$, which is a function on $X$.

The aspect of K\"ahler geometry which is our concern here is the existence of particular metrics, in a given K\"ahler class, which correspond to the solutions of partial differential equations for the potential function $\phi$. The four equations we have in mind are as follows.
\begin{itemize}\item {\it K\"ahler-Einstein metrics} with $\rho=\lambda \omega$ for $\lambda=1,-1$ or $0$. These can only occur when $c_{1}(X)$ (a topological  invariant of the complex structure) is the corresponding multiple of the K\"ahler class.
\item {\it Constant scalar curvature K\"ahler} (CSCK) metrics with $S={\rm constant}$. The constant is determined by the topological data, $c_{1}(X),  [\omega]$. When $c_{1}(X)$ is a multiple of $[\omega]$ an integral identity shows that a CSCK metric is K\"ahler-Einstein.
\item {\it Extremal metrics}. The satisfy the condition that the gradient of $S$, a vector field on the manifold, is holomorphic.  Obviously CSCK metrics (where the gradient is zero) are extremal and in general the holomorphic vector field is determined by the topological data. 
\item {\it K\"ahler-Ricci solitons}.  These satisfy the condition that
$$  \rho-\lambda\omega= L_{v}\omega$$
where $v$ is a holomorphic vector field and $L_{v}$ denotes the Lie derivative. Obviously K\"ahler-Einstein metrics satisfy this condition (with $v=0$) and, as before,   $v$ is determined by the topological data.  \end{itemize}

Extremal metrics and K\"ahler-Ricci solitons arise naturally in the study of associated parabolic equations. The {\it K\"ahler-Ricci flow} is the evolution equation
$$    \frac{\partial \omega}{\partial t}= \lambda \omega- \rho. $$
A K\"ahler-Ricci soliton represents a fixed point of this flow modulo the holomorphic diffeomorphisms---the flow starting from a Kahler-Ricci soliton evolves as a family of geometrically equivalent metrics. The {\it Calabi flow} is the evolution equation
$$  \frac{\partial \omega}{\partial t} = i\dbd S, $$
and extremal metrics appear in the same way.

\subsection{The YTD conjecture}

The existence of K\"ahler-Einstein metrics in the cases when $\lambda=-1$ and $\lambda=0$  was established in the 1970's in breakthrough work of Yau \cite{kn:Y} confirming conjectures of Calabi from the 1950's \cite{kn:Cal1},\cite{kn:Cal2}. (The case when $\lambda=-1$ was treated independently by Aubin.) The positive case, when $\lambda=1$ (which arises when $X$ is a Fano manifold) is more subtle and, beginning with a theorem of Matsushima---also from the 1950's---various obstructions were found to the existence of a K\"ahler-Einstein metric. In the 1980's Yau suggested that the existence of these metrics should be equivalent to some algebro-geometric notion of \lq\lq stability''. One motivation for this idea came from the theory of {\it Hermitian Yang-Mills connections}.   Let $E$ be a holomorphic vector bundle over a K\"ahler manifold $(X,\omega_{0})$ and suppose for simplicity that $c_{1}(E)=0$ and that $E$ is indecomposable (cannot be written as a direct sum). A Hermitian Yang-Mills connection on $E$ is a unitary connection, compatible with the holomorphic structure, such that the inner product of the curvature with the K\"ahler form vanishes. The bundle $E$ is called {\it stable} if  any proper subsheaf ${\cal S}\subset {\cal O}(E)$ has $c_{1}({\cal S})<0$. This is a purely complex-geometric notion (or algebro-geometric, in the case when $X$ is a projective  variety), which arose in the study of moduli problems for holomorphic bundles.  The main result is that stability is the necessary and sufficient condition for the existence of a Hermitian Yang-Mills connection. This was proved by Uhlenbeck and Yau \cite{kn:UY} in 1986 (with an independent treatment by Donaldson in the algebraic case), confirming conjectures made by Kobayashi and Hitchin.  Many extensions of the result have been found,  for example to coupled equations   with connections and additional \lq\lq fields''. Yau's proposal was that the K\"ahler-Einstein question could be put in a similar framework.
 
 Considerable progress towards confirming Yau's proposal was made by Tian during the 1990's, In particular, Tian gave a definition of the notion of stability for Fano manifolds, called {\it K-stability}, which turns out to be the correct one \cite{kn:Ti}. Around the turn of the century, Donaldson suggested that the whole conjectural picture could be extended to the question of existence of constant scalar curvature metrics in a given rational K\"ahler class on an projective variety, with a suitable extension of the definition of K-stability \cite{kn:D-2}, \cite{kn:D-1}.  Later, Sz\'ekelyhidi proposed further extensions to the existence of extremal K\"ahler metrics \cite{kn:Gabor0}. 

We now recall the definition of K-stability. We restrict attention to the case when the K\"ahler class is the first Chern class of a positive line bundle $L$, so in algebro-geometric terms we have a polarised manifold $(X,L)$. The definition is founded on the notions of {\it test configurations} and {\it Futaki invariants}. A test configuration  for $(X,L)$ is a normal variety  ${\cal X}$ which forms a flat family $\pi: {\cal X}\rightarrow \bC$ with a line bundle ${\cal L}\rightarrow {\cal X}$ such that:
\begin{itemize}
\item ${\cal L}$ is ample on the fibres of $\pi$;
\item there is a $\bC^{*}$-action on $({\cal X},{\cal L})$ covering the standard action on $\bC$;
\item the fibre $\pi^{-1}(1)$ is isomorphic to $X$ and ${\cal L}\vert_{\pi^{-1}(1)}=L^{r}$ for some $r$. 
\end{itemize}
Of course the $\bC^{*}$-action means that all non-zero fibres are isomorphic to $X$ and the essential point is that a test configuration gives a {\it degeneration} of $X$ to a typically-different object: the scheme-theoretic central fibre $X_{0}=\pi^{-1}(0)$. For example we might degenerate a smooth conic in the plane to a pair of lines or a \lq\lq double line''. 

Given such a test configuration, the  central fibre $X_{0}=\pi^{-1}(0)$, with the restriction ${\cal L}_{0}$, is a polarised scheme with $\bC^{*}$-action. In such a situation there is a numerical invariant, the Futaki invariant, which can be defined in various ways. One is through the induced action on a \lq\lq CM line'' associated to $(X_{0}, {\cal L}_{0})$; another is through the asymptotics as $k\rightarrow \infty$ of the induced action on the vector spaces $H^{0}(X_{0}, {\cal L}_{0}^{k})$. The upshot is that we get a numerical invariant ${\rm Fut}({\cal X})$ of a test configuration.  We say that $(X,L)$ is K-stable if ${\rm Fut}({\cal X})\geq 0$ for all  test configurations  and equality occurs only if ${\cal X}$ is a product $X\times \bC$, with a possibly non-trivial $\bC^{*}$ action on the $X$ factor. The \lq\lq standard conjecture'' in the field,  which has come to be called the Yau-Tian-Donaldaon (YTD) conjecture, is that there is a constant scalar curvature metric on $X$ in the class $c_{1}(L)$ if and only if $(X,L)$ is K-stable. An inkling of the connection between the algebro-geometric and differential geometric sides is given by considering the case when the central fibre $X_{0}$ is smooth. Choose a K\"ahler form in the class $c_{1}({\cal L}_{0})$ on $X_{0}$ invariant under the action of the subgroup $S^{1}\subset \bC^{*}$. The action  is  generated by a Hamiltonian function $H$ on $X_{0}$, with respect to this K\"ahler form, regarded as a symplectic structure. One then has a formula
  $$    {\rm Fut}({\cal X})= \int_{X_{0}} ( S-\hat{S})\  H \ d\mu, $$
  where $S$ is the scalar curvature of $X_{0}$ and $\hat{S}$ is the average value of $S$.

 This YTD conjecture has close parallels with the Kobayashi-Hitchin conjecture for Hermitian-Yang-Mills connections discussed above and in fact the definition of stability for holomorphic bundles can be put in a similar shape, involving suitable degenerations of the vector bundle. 

\

\

The YTD conjecture has been confirmed in certain cases, in particular for toric surfaces \cite{kn:D0} and for Fano manifolds, as we will discuss below. However there are reasons to think that in general the conjecture may not be exactly true as formulated and that the correct conjecture should involve a slightly different definition of stability. One such notion has been proposed by Sz\'ekelyhidi \cite{kn:Gabor1}. This involves an algebraic generalisation of the notion of a test configuration having the form of a filtration of the co-ordinate ring 
$$   R_{X,L}= \bigoplus_{k} H^{0}(X,L^{k}). $$
The set of filtrations can be thought of roughly as a \lq\lq completion'' of the set of test configurations. Sz\'ekelyhidi obtained a slightly stronger notion called $\hat{K}$-stability, using an extension of the definition of the Futaki invariant to filtrations.

Another important variant is \lq\lq uniform K-stability'' which also goes back to work of Sz\'ekeleyhidi, with later developments by Dervan \cite{kn:Dervan} and Boucksom, Hisamoto, Jonsson \cite{kn:BHJ}. This notion depends upon a choice of a \lq\lq norm'' $\Vert {\cal X}\Vert$ on test configurations. For example in the case of a smooth central fibre considered above this could be given by an $L^{p}$ norm of the Hamiltonian function $H$,  normalised to have integral zero: the definition in the general case is more technical.
 A pair $(X,L)$ with no non-trivial automorphisms is called uniformly K-stable if there is some $\epsilon>0$ such that for all test configurations ${\cal X}$
$$   {\rm Fut}({\cal X})\geq \epsilon \Vert {\cal X}\Vert. $$
(The definition needs to be modified slightly in the case when ${\rm Aut}(X,L)$ is non-trivial.)

The relation between these different notions of stability can be illustrated by considering the case of toric pairs $(X,L)$. Such  pairs corresponds to certain closed polytopes $P\subset \bR^{n}$, which come with a distinguished measure $d\sigma$ on the boundary. Then we can define a linear functional on functions on $P$,
$$   L_{P}(f)= \int_{\partial P} f d\sigma - A \int_{P} f d\mu, $$
where the constant $A$ is fixed by requiring that $L_{P}$ vanishes on constant functions. We say that a function $f$ is \lq\lq normalised'' if 
$$  \int_{P} f\  d\mu=0  ,  \int_{P} f \ x_{i} \ d\mu =0 $$
for each co-ordinate function $x_{i}$. We say that $f$ is a {\it convex rational piecewise-linear}  function if
$  f=\max \{ \lambda_{r} \}$,
where $\{\lambda_{r}\}$ is a finite collection of affine-linear functions with rational coefficients. Write $QPL$ for the set of such functions $f$.
Then we may consider three properties of the polytope $P$. 
\begin{enumerate}\item For all non-zero normalised $f\in QPL$  we have $L_{P}(f)>0$.
\item For all non-zero normalised continuous convex functions $f$ on $P$ we have $L_{P}(f)>0$.
\item There is an $\epsilon>0$ such for all normalised $f\in QPL$ we have
  $$   L_{P}(f) \geq \epsilon \int_{P} \vert f\vert d\mu. $$
 \end{enumerate}

These three properties correspond respectively to K-stability,  $\hat{K}$-stability and uniform K-stability of the pair $(X,L)$ determined by $P$. (Notice that in the last item it does not matter whether we consider rational piecewise linear or general continuous convex functions, since the former are a dense subset of the latter.) 

\subsection{Fano manifolds and K\"ahler-Einstein metrics}
 The YTD conjecture in the case of Fano manifolds and K\"ahler-Einstein metrics (the case originally considered by Yau) was confirmed in \cite{kn:CDS}. Since this work was discussed by Sz\'ekelyhidi in his contribution to the 2014 ICM \cite{kn:Gabor2} we will not say much about it here. By now, four different proofs of (essentially) this result have appeared; the references for the other three being \cite{kn:DaS}, \cite{kn:CSW}, \cite{kn:BBJ}. Three of these proofs follow somewhat similar strategies. The approach is to take some deformation or continuity procedure and show that either it produces a K\"ahler-Einstein metric on our Fano manifold $X$ or it produces (roughly speaking) the central fibre $X_{0}$ in a test configuration which contradicts K-stability. These three proofs depend heavily on the Cheeger-Colding convergence theory for manifolds with Ricci curvature bounds which we will come back to in Section 3.4 below.  In the first proof, by Chen, Donaldson and Sun,  the deformation process involves  K\"ahler-Einstein metrics with cone singularities along a divisor in $X$ and increasing the cone angle from some small initial value. If we achieve cone angle $2\pi$ we have a smooth metric. In the second proof, Datar and Sz\'ekelyhidi used the Aubin-Yau continuity method with the $1$-parameter family of equations  
$$  {\rm Ric}(\omega_{s})= (1-s)\alpha + s \omega_{s},  $$
where $\alpha$ is a fixed positive $(1,1)$-form. In the third proof, Chen and Wang \cite{kn:CW} study the K\"ahler-Ricci flow $\omega_{t}$ starting with some initial metric on $X$. The flow  exists for all positive time and Chen and Wang show that for any sequence $t_{i}\rightarrow \infty$ there is a subsequence $t'_{i}$ such that the Riemannian manifolds  $(X,\omega_{t'_{i}})$ converge geometrically to a K\"ahler-Ricci soliton metric on some possibly singular space $Z$,  which is a complex algebraic Q-Fano variety. In \cite{kn:CSW} Chen, Sun and Wang showed further that this limit is unique,  so in fact the K\"ahler-Ricci flow $(X,\omega_{t})$ converges in this sense as $t\rightarrow\infty$. The case when $Z=X$ gives the desired K\"ahler-Ricci soliton metric on $X$ (which is K\"ahler-Einstein if $X$ is K-stable), and if $Z$ is not equal to $X$ then Chen, Sun and Wang show how to produce a destabilising test configuration from $Z$. 

The information which these proofs using convergence theory provide about the precise way in which the possible  failure of the PDE construction strategies relates to K-stability has independent interest. There are parallel questions in the case of the Hermitian-Yang-Mills equations which have been fully settled relatively recently. Jacob\cite{kn:J} and Sibley \cite{kn:Sibley} study the Hermitian Yang-Mills flow (analogous to the K\"ahler-Ricci flow) and show that the limit of this produces the algebro-geometric \lq\lq Harder-Narasimhan stratification'' of a holomorphic bundle, thus refining the basic existence result for Hermitian-Yang Mills connections on stable bundles. 

\

The fourth proof of the YTD conjecture for Fano manifolds, by Berman, Boucksom and Jonsson, uses very different techniques. It relies  on foundations from algebraic geometry and pluripotential theory, rather than from Riemannian geometry,  and they exploit the variational point of view which we discuss in 2.4 below.  The result proved is slightly weaker: they show when $X$ has finite automorphism group the existence of a K\"ahler-Einstein metric is equivalent to uniform K-stability.

\

Over the past four years the interaction between \lq\lq K\"ahler-Einstein geometry''
 and algebraic geometry has blossomed, with many related developments. One important theme is given by various notions related to \lq\lq volume''. A beautiful example  is a result of Fujita \cite{kn:Fuj}. In differential geometric language this states that if a compact complex $n$-manifold $X$ has a K\"ahler-Einstein metric with  ${\rm Ricci}=\omega$ then the volume of the manifold is at most $(2\pi)^{n} (n+1)^{n}$, with equality if and only if $X=\bC\bP^{n}$ (with its usual Fubini-Study metric, suitably scaled). This can be set alongside the standard Bishop comparison result, which states that among all Riemannian manifolds with ${\rm Ricci}= g$ (in fact with ${\rm Ricci}\geq g$) the round sphere (of the right scale) maximises the volume. So Fujita's result is a K\"ahler analogue.   But while the statement is differential geometric, Fujita's proof is  entirely algebro-geometric, using the equivalence with $K$-stability. In the other direction,  differential geometric results on K\"ahler-Einstein metrics have had important consequences for the construction of moduli spaces of Fano manifolds in algebraic geometry, giving information which cannot at present be obtained by purely algebro-geometric methods. We will not go into this important topic further here but refer to the contribution of Song Sun to these ICM Proceedings \cite{kn:Song}.

Another line of development involves singularities. An important part of the Cheeger-Colding convergence theory is that the relevant limit spaces have metric {\it tangent cones} at each point. In a setting (such as the limits of the Ricci flow discussed above, or in considering compactified moduli spaces) where one has such a limit space $Z$ which is a singular complex algebraic variety one can ask how this metric theory is related to the algebraic geometry of the singularities. This question was studied by the author and Sun \cite{kn:DS}, who showed that the metric tangent cone at a point $p\in Z$ is unique and can be obtained from a valuation $\nu_{p}$ on the local ring of $Z$ at $p$. This valuation records the order of vanishing of a function,  with respect to the metric on $Z$:
$$  \nu_{p} (f) = \lim_{r\rightarrow 0} \frac{\log \max_{B_{r}} \vert f\vert}{\log r}, $$
where $B_{r}$ is the $r$-ball in $Z$ centred at $p$. The work of Donaldson and Sun left open the question of whether the valuation can be determined entirely algebro-geometrically. This was settled by Hein and Sun \cite{kn:HS} in the case of singular Calabi-Yau varieties with certain singularities, including ordinary double points (which is the most important case for applications). Much progress in the general case has been made by Li, Liu and Xu, who show that the valuation can be characterised as minimising a certain algebro-geometrically defined normalised volume \cite{kn:Li}, \cite{kn:LiLiu},\cite{kn:LiuXu}, \cite{kn:LiXu2} Again we refer to \cite{kn:Song} for a more detailed account of these developments.

\

The tangent cones that appear in this context are, differential-geometrically, {\it Sasaki-Einstein cones} (possibly with singularities). The algebro-geometric setting is to consider an action of a torus $T=\left(\bC^{*}\right)^{r}$ on $\bC^{N}$ and  a T-invariant affine variety $W\subset \bC^{N}$ with trivial canonical bundle.  The general existence question in this setting asks if there is a Sasaki-Einstein metric on  $W$ with radial vector field generated by a vector in the Lie algebra of $T$.    One case is when $r=1$ and we have the standard action of $\bC^{*}$,  so $W$ is a cone in the algebro-geometric sense over a  projective variety and the existence problem reduces to the existence of a K\"ahler-Einstein metric on this projective variety.  The analogue of the YTD conjecture in this Sasaki-Einstein setting (for varieties $W$ with an isolated singularity at $0$) was established by Collins and Sz\'ekelyhidi \cite{kn:CS}. 

\

A glaring problem in this field is that it is usually very hard to decide if a manifold is K-stable. On the face of it, the definition requires one to check all possible test configurations. (This is in contrast to the situation for bundles, where the criterion for stability is very explicit.) We mention three lines of recent progress in this direction.

\

 \begin{itemize}\item The proofs of the YTD conjecture show that in the Fano case it suffices to check test configurations with normal central fibres (as in Tian's original proposal). This was also proved algebro-geometrically a little earlier by Li and Xu \cite{kn:LiXu}.
 \item For explicit examples of interest such as cubic 3-folds in $\bP^{4}$ \cite{kn:LiuXu}, and intersections of two quadrics \cite{kn:SSun}, K-stability is completely understood (including for the singular varieties that are added to compactify moduli spaces).
\item Various cases of manifolds with large symmetry groups have been analysed \cite{kn:IS}, \cite{kn:Delcroix2}.
\item Fujita and Odaka \cite{kn:FujOdaka} defined a numerical invariant $\delta(X)$ of a Fano manifold and conjectured that uniform K-stability is equivalent to $\delta(X)>1$, which was later confirmed by Blum and Jonsson \cite{kn:BlumJ}. This invariant $\delta$ is obtained from the theory of the log canonical threshold (lct) (which in this context is a numerical invariant measuring how singular a divisor is).   In fact $\delta(X) =\limsup_{k\rightarrow\infty} \delta_{k}(X)$ with
$$ \delta_{k}= \inf_{D} {\rm lct}(X,D).$$
Here, for each $k$, the divisor $D$ ranges over sums $\sum_{i} D_{i}$, of the divisors associated to a basis $s_{i}$ of $H^{0}_{k}= H^{0}(X, K_{X}^{-k})$ and with a normalising factor $\left(k \ {\rm dim} H^{0}_{k}\right)^{-1}$.  Fujita and Odaka also relate their invariant $\delta$ to another definition of \lq\lq Gibbs stability'' proposed by Berman.
\end{itemize}
\subsection{Geometry in the space of K\"ahler metrics}

The existence questions in K\"ahler geometry which we have been discussing have appealing formulations in terms of infinite-dimensional geometry. Given the K\"ahler class $[\omega_{0}]$ we write ${\cal H}$ for the space of K\"ahler potentials $$  {\cal H}= \{\phi: \omega_{0}+ i\dbd \phi>0\}. $$
(When the K\"ahler class is integral this has a more invariant description in terms of metrics on the corresponding ample line bundle.) We define a Riemannian metric, making ${\cal H}$ into an infinite-dimensional Riemannian manifold
$$  \Vert \delta\phi\Vert^{2}_{\phi} = \int_{X} (\delta\phi)^{2} d\mu_{\phi}, $$
where $d\mu_{\phi}$ is the usual volume form $\omega_{\phi}^{n}/n!$. This metric was discovered by Mabuchi \cite{kn:Mab1} and rediscovered later by Semmes and Donaldson. It has the remarkable property that ${\cal H}$ becomes an infinite-dimensional symmetric space of non-positive curvature, formally the dual of the group of Hamiltonian symplectomorphisms of $(X,\omega_{0})$ \cite{kn:D-1}. A geodesic segment in ${\cal H}$ corresponds to an $S^{1}$-invariant degenerate $(1,1)$ form $\Omega$ on $(0,1)\times S^{1}\times X$ with $\Omega^{n+1}=0$, which in turn corresponds to a solution of a certain homogeneous complex Monge-Amp\`ere equation. 

Various functionals on ${\cal H}$ play an important role in the theory. Typically these are defined through their first variations (that is, what is defined is a closed $1$-form on ${\cal H}$). The {\it Mabuchi functional} ${\cal F}$ is defined by
$$  \delta {\cal F}= \int_{X} (\delta \phi) (S-\hat{S}) \ d\mu_{\phi},  $$
where $S$ is the scalar curvature and $\hat{S}$ is the average \cite{kn:Mab0} By construction, the Euler-Lagrange equation associated to the Mabuchi functional is the constant scalar curvature condition. The Mabuchi functional is convex along geodesics in ${\cal H}$; so in a formal way the search for a constant scalar curvature metric is the search for a critical point (in fact a global minimum) of a convex function on a space of negative curvature. This formal picture can be fitted into a general framework involving \lq\lq stability'' and moment maps---see the exposition in \cite{kn:D2},  for example. 

\

The first steps in making this formal picture of real use in K\"ahler geometry were taken by Chen \cite{kn:XC}. He showed that any two points in ${\cal H}$ can be joined by a geodesic segment, provided that the definitions were extended to allow non-smooth potentials, roughly speaking of class $C^{1,1}$. Calabi and Chen showed in \cite{kn:CC} that there is a genuine induced metric on ${\cal H}$, and that it is a space of negative curvature in a generalised sense of triangle comparison. Chen also found an important formula for the Mabuchi functional.  Let $V_{\phi}$ be the volume element of $\omega_{\phi}$ in terms of the fixed reference
metric $\omega_{0}$:
$$  V_{\phi} \omega_{0}^{n}= \omega_{\phi}^{n}$$ and let $\rho_{0}$ be the Ricci form of $\omega_{0}$. Then Chen's formula is

$$  {\cal F}(\phi)= \int_{X} V_{\phi}\log V_{\phi}\  d\mu_{0} + J_{\rho_{0}}(\phi), $$ where the \lq\lq J-functional'' is defined by
$$  \delta J_{\rho_{0}} = \int_{X} (\delta \phi)\  \rho_{0} \wedge \omega_{\phi}^{n-1}. $$ One important consequence of this formula is that it gives a definition of ${\cal F}$ for $C^{1,1}$ potentials. 

The connection with K-stability and algebraic geometry comes from the fact that under certain conditions a test configuration gives a geodesic ray in ${\cal H}$ and the asymptotic behaviour of the Mabuchi functional along the ray can be related to the Futaki invariant of the test configuration. Statements along these lines, under various technical hypotheses,  have been proved by Phong and Sturm \cite{kn:PS} and many other authors.

\

There have been important developments in this area over the past few years. One is due to Berman and Berndtsson \cite{kn:BB},  who showed that ${\cal F}$ is convex along generalised geodesics.  This leads to a short and conceptual proof of various foundational results, such as the uniqueness of constant scalar curvature metrics. (All of this discussion can be modified to apply to extremal metrics.) Another is due to Darvas \cite{kn:Darvas}, who showed that a metric completion of ${\cal H}$ has a good analytical meaning. In fact it turns out to be best to start with a Finsler metric on ${\cal H}$, defined by the $L^{1}$-norm rather than the $L^{2}$-norm. Then Darvas showed that the completion can be identified with a space ${\cal E}^{1}$ of \lq\lq finite-energy'' currents which was introduced previously in \cite{kn:BBEGZ} and which is important in pluripotential theory.

These ideas form some of the background for the version of the YTD conjecture proved by Berman {\it et al}. In fact they work with both the Mabuchi functional and the {\it Ding functional}, which has a similar character but is only defined
in the K\"ahler-Einstein context. Very roughly, the strategy is as follows.
\begin{itemize}
\item Existence of the minimum follows from properness of the Ding functional;
\item If the functional is not proper there is a geodesic ray on which the functional is bounded above;
\item This geodesic ray can be approximated by test configurations and the boundedness gives a contradiction to uniform K-stability.
\end{itemize}

The same kind of strategy can be applied beyond the K\"ahler-Einstein case. For example Darvas and Rubinstein show that (for any K\"ahler class), properness of the Mabuchi functional (with respect to the $L^{1}$-distance) implies the existence of \lq\lq weak'' minimisers, in the completion ${\cal E}^{1}$ \cite{kn:DR}. This brings to the fore the question of regularity of the weak minimisers. In a similar vein, Streets \cite{kn:Streets} defined a generalised \lq\lq minimising movement'' version of the Calabi flow (exploiting on the favourable metric geometry of ${\cal H}$),  and proved that this satisfies long time existence. So there  are parallel question about the regularity of this Streets flow.  

 Very recent work of Chen and Cheng makes important progress in this variational and PDE framework \cite{kn:ChenCheng}. In connecting with these developments with algebraic geometry, one central question is whether uniform K-stability implies properness of the Mabuchi functional.

\section{Exceptional holonomy}
\subsection{Background}
Let $(M,g)$ be a Riemannian manifold of dimension $m$ and $p\in M$ a basepoint. Parallel transport around a loop based at $p$ defines an orthogonal transformation of the tangent space at $p$. The  set of all these transformations arising from contractible loops is a closed, connected, subgroup of the orthogonal group called the {\it holonomy group} of $(M,g)$. This can be viewed, up to conjugacy, as a subgroup $G$ of $SO(m)$.  Work of Berger from the 1950's, with refinements by later authors, gives a complete classification of all groups that can arise. Leaving aside reducible subgroups (where the manifold has a local product structure) and symmetric spaces (which are completely classified) there are five standard families and two exceptional cases. The standard families fall into three pairs.
\begin{itemize}
\item {\bf Real}: $G=SO(m)$. 
\item {\bf Complex}: $m$ is even and $G=U(m/2)$ or $G=SU(m/2)$.
\item {\bf Quaternionic}: $m$ is divisible by $4$ and $G=Sp(m/4)$ or $G=Sp(m/4). Sp(1)$
\end{itemize}

The first case is that of a generic Riemannian metric. In the second we make the standard identification $\bR^{m}= \bC^{m/2}$. In the third we make the standard identification $\bR^{m}= \bH^{m/4}$ where $\bH$ denotes the quaternions and the group $Sp(n).Sp(1)$ is given by left multiplication by quaternionic matrices in $Sp(n)$ and right multiplication by unit quaternions $Sp(1)$.

The exceptional cases occur in dimensions $7$ and $8$ and correspond to the exceptional Lie group $G_{2}\subset SO(7)$ and the subgroup ${\rm Spin}(7)\subset SO(8)$ given by the spin representation in dimension $7$. Just as the standard families are  built on the algebra of the three fields $\bR,\bC,\bH$, the exceptional cases are built on \lq\lq exceptional'' algebraic phenomena, associated to the Octonions, triality, etc. The basic definitions can be approached in various   ways: a convenient approach for differential geometers emphasises exterior algebra.

Let $V$ be a real vector space of dimension $m$. For $m\geq 9$ and $2<p<m-2$ the dimension of the exterior power $\Lambda^{p}V^{*}$ exceeds that of the general linear group $GL(V)$. But in low dimensions there are a few cases where $GL(V)$ acts with an open orbit. In particular this happens when $m=7$ and $p=3$.  For an oriented 7-dimensional $V$ and $\phi\in \Lambda^{3}V^{*}$ we consider the quadratic form on $V$
\begin{equation}  v\mapsto i_{v}(\phi)\wedge i_{v}(\phi)\wedge \phi. \end{equation}
This takes values in the oriented line $\Lambda^{7}V^{*}$, so it makes sense to say that the form is positive definite. In such a case we say that $\phi$ is a positive $3$-form. The basic facts are:
\begin{itemize}\item the positive $3$-forms make up a single orbit for the action of $GL^{+}(V)$;
\item the stabiliser in $GL^{+}(V)$ of a positive $3$-form is a compact Lie group isomorphic to the exceptional Lie group $G_{2}$.
\end{itemize}
With regard to the second item; a positive $3$-form $\phi$ determines a Euclidean form $g_{\phi}$ on $V$, so the stabiliser is a subgroup of the orthogonal group for this Euclidean structure. Indeed the quadratic form (1) gives a conformal class of Euclidean structures, so the definition of $g_{\phi}$ is just a matter of fixing the scale. This can be achieved by requiring that $\vert \phi\vert^{2}=7$ in the norm induced by $g_{\phi}$. Let $*_{\phi}$ be the $*$-operator determined by $g_{\phi}$. Then we have a 4-form $*_{\phi}\phi$ which is also preserved by the stabiliser of $\phi$.

We now go to dimension 8 by considering the vector space $W=V\oplus \bR$ with co-ordinate $t$ in the $\bR$ factor. Given a positive $3$-form $\phi$ as above we consider the $4$-form
$$      \Omega= \phi\wedge dt + *_{\phi}\phi \in \Lambda^{4}(W^{*}). $$
The basic fact is that the stabiliser of $\Omega$ in $GL^{+}(W)$ is a compact Lie group isomorphic to ${\rm Spin}(7)$ (the double cover of $SO(7)$). We
 write ${\cal A}$ for the orbit under $GL^{+}(W)$ of $\Omega$ in $\Lambda^{4}W^{*}$, so ${\cal A}$ has dimension $64-21=43$. The forms in ${\cal A}$ are called {\it admissible}.

Moving  to differential geometry, we have a notion of a positive 3-form on an oriented 7-manifold and an admissible 4-form on an oriented 8-manifold. The link with exceptional holonomy is provided by a result of Fern\'andez and Gray \cite{kn:FGray}.  To state this precisely it is useful to make a small change of viewpoint 
and to talk about a \lq\lq torsion-free $G$-structure'', which is to say
a reduction of the structure group of the tangent bundle of the manifold
to a group $G\subset SO(m)$ and a  torsion-free $G$-connection on the tangent
bundle. This  gives a Riemannian metric with holonomy contained in $G$  and the notion takes care of technical complications when considering holonomy around non-contractible loops and when the holonomy is strictly smaller than $G$. The result of Fern\'andez and Gray is:

  \begin{itemize} \item Giving a torsion-free $G_{2}$-structure on an oriented $7$-manifold  is equivalent to giving a positive 3-form $\phi$ such that $d\phi=0$ and $d*_{\phi} \phi =0$.
\item Giving a torsion free ${\rm Spin}(7)$-structure on an oriented $8$-manifold is equivalent to giving a  admissible 4-form $\Omega$ with $d\Omega=0$.
\end{itemize}

There are many other points of view on these two exceptional holonomy groups. For example, they also can be characterised by the existence of a covariant constant spinor field. A very important feature is that the Ricci curvature is identically zero in both cases. 

\

For the rest of this article we will concentrate on the 7-dimensional case. For euphony here we will  call  a torsion-free $G_{2}$-structure simply a $G_{2}$-structure.

\subsection{Examples of manifolds with holonomy $G_{2}$}

While the Berger classification goes back to the 1950's, the modern developments in this area  begin with work of Bryant in the mid 1980's \cite{kn:Br1}, giving the first local examples of metrics with exceptional holonomy--see also Bryant's lecture in ICM 1986 \cite{kn:Br2}. 
The global theory took off about 10 years later when Joyce established the existence of compact examples \cite{kn:Joyce1}, \cite{kn:Joyce2}.  Joyce's approach was to start with a singular quotient space such as $T/\Gamma$ where $T$ is a flat torus and $\Gamma$ is a finite group and then remove the singularities by a \lq\lq gluing'' construction. This work was described in Joyce's 1998 ICM lecture \cite{kn:Joyce3} and in a monograph \cite{kn:Joyce4} which gives a comprehensive treatment of the field up to the turn of the century. There is at present essentially only one other source of examples of compact manifolds with $G_{2}$-holonomy, given by the \lq\lq twisted connected sum'' construction introduced by Kovalev in \cite{kn:Kov1}.

  We  give a brief outline of these Joyce and Kovalev constructions. One might think of these as being based on two models for a positive $3$-form, corresponding to the two subgroups: $SO(4)\subset G_{2}$ and $SU(3)\subset G_{2}$.  For the first model  we take $\bR^{4}$, regarded as an oriented 4-dimensional Euclidean vector space, and choose a standard basis $\omega_{1}, \omega_{2}, \omega_{3}$ for the space of self-dual $2$-forms $\Lambda^{2}_{+}$ on $\bR^{4}$. Then we take $\bR^{3}$ with co-ordinates $t_{1}, t_{2}, t_{3}$ and write down the positive $3$-form
  \begin{equation}  \phi= dt_{1} dt_{2} dt_{3} - \sum \omega_{i} dt_{i} \end{equation}
  on the $7$-dimensional vector space $\bR^{4}\oplus \bR^{3}$. In more invariant terms, there is a canonical positive 3-form on $\bR^{4}\oplus \Lambda^{2}_{+}$, preserved by the action of $SO(4)$. Now let $Q$ be a {\it hyperk\"ahler 4-manifold}, which is the same as saying that there is a orthonormal frame of closed self-dual forms $\omega_{i}$ on $Q$. The same formula (2) defines a closed and co-closed positive $3$-form on on $Q\times \bR^{3}$ or $Q \times T^{3}$, which of course corresponds to a metric with holonomy strictly contained in $G_{2}$. 
The simplest case of Joyce's construction arises from a quotient $T^{7}/\Gamma$ where the singular set is a disjoint union of flat 3-tori and a neighbourhood of each component is modelled on $T^{3}\times \bR^{4}/\pm{1} $.
The Eguchi-Hanson manifold $Q$ is a non-compact hyperk\"aher 4-manifold asymptotic to $\bR^{4}/\pm 1$. Joyce removes neighbourhoods of each of the components of the singular set and replaces them by corresponding neighbourhoods in $T^{3}\times Q_{\epsilon}$, where $Q_{\epsilon}$ is the same manifold $Q$ with the metric scaled by a factor $\epsilon$. He defines a positive 3-form---the \lq\lq approximate solution''--- on the resulting 7-manifold by gluing the structures on each piece and then shows, using a PDE and implicit function theorem argument, that for small $\epsilon$ this approximate solution can be deformed slightly to obtain a  $G_{2}$-structure. 

\

For the second model we start with $\bC^{3}$, regarded as having symmetry group $SU(3)$, so we have the standard Hermitian metric with $(1,1)$ form $\omega$ and a standard complex 3-form $\Theta$. Thus, in the usual co-ordinates $z_{a}=x_{a}+ i y_{a}$:
$$  \omega= \sum dx_{a} dy_{a}\ \ \ , \ \ \ \Theta= dz_{1}dz_{2}dz_{3}. $$
We write down a positive 3-form 
\begin{equation}   \omega\wedge dt + {\rm Re}(\Theta), \end{equation}
on the 7-dimensional real vector space $\bC^{3}\oplus \bR$, where $t$ is the coordinate on the $\bR$-factor. Let $Z$ be a 6-manifold with holonomy $SU(3)$. This is the same as saying that $Z$ is a three complex-dimensional K\"ahler manifold with K\"ahler form $\omega$ and a holomorphic $3$-form $\Theta$ of constant norm: at each point we can choose co-ordinates so that the forms match up with the standard models above. (This is also called a {\it Calabi-Yau} structure on $Z$. Up to coverings, it is equivalent to saying that $Z$ has a K\"ahler-Einstein metric with zero Ricci curvature: the case $\lambda=0$ in the language of the first part of this article.) The same formula (3) defines a $G_{2}$ structure on the product $Z\times \bR$ or $Z\times S^{1}$.  

Kovalev's construction depends on a supply of suitable examples of Calabi-Yau manifolds $Z$,  and these come from a  variant of the K\"ahler theory discussed in the first part of this article. Kovalev takes a  Fano 3-fold $W$ with a smooth anticanonical divisor $D\subset W$. By standard algebraic geometry theory, $D$ is a complex $K3$ surface. (An example is $W=\bC\bP^{3}$ and $D$ a smooth quartic surface.) Also by standard theory, we can choose  a curve $C\subset D$ so that the proper transform $\tilde{D}$ of $D$ in the blow-up $\tilde{W}$ of $W$ along $C$ has trivial normal bundle. (Of course, the proper transform $\tilde{D}$ is isomorphic to $D$.) Now let $Z$ be the complement of $\tilde{D}$ in $\tilde{W}$. A result of Tian and Yau, in the general vein of the theory discussed in the first part of this article, gives the existence of a complete Calabi-Yau structure on Z. Kovalev showed (with later clarification  by Haskins, Hein and Nordstr\"om \cite{kn:HHN}) that this Tian-Yau metric is {\it asymptotically cylindrical}--- it is asymptotic to a cylinder $\tilde{D}\times S^{1}\times \bR$ where the K3 surface $\tilde{D}=D$ is endowed with the Ricci-flat K\"ahler-Einstein metric given by Yau's theorem. Thus, from the preceding discussion, we get an asymptotically cylindrical 7-manifold $Z\times S^{1}$ with holonomy strictly contained in $G_{2}$ and with asymptotic cross-section 
\begin{equation}\tilde{D}\times S^{1}\times S^{1}. \end{equation} 

Next, Kovalev considers two manifolds $Z_{1}, Z_{2}$ constructed by the above procedure. These have asymptotic cross-sections $\tilde{D}_{i}\times S^{1}\times S^{1}$. He cuts off the ends of these manifolds, taking a large parameter $L$ which will be roughly the diameter of the remaining pieces,  and constructs a compact $7$-manifold $M$ by gluing these pieces together. The \lq\lq twist'' in the construction is that in making this gluing he interchanges the two $S^{1}$ factors in the cross-sections (4) , so there is no global $S^{1}$ action on  $M$. To match up (asymptotically, for large $L$) the positive 3-forms on the two pieces the $K3$ surfaces $\tilde{D}_{i}$ must satisfy a certain matching condition involving \lq\lq hyperk\"ahler rotation'' of the complex structures.  This matching problem can be studied through the highly-developed  Torelli theory on $K3$ surfaces, and the upshot was that Kovalev was able to find many examples satisfying the condition. Gluing the positive three forms on the two pieces in $M$ gives an approximate solution and Kovalev showed that this could be deformed to a genuine $G_{2}$-structure on $M$, for large values of the parameter $L$.

\

The Joyce and Kovalev constructions have similarities. They both start with building blocks constructed  using better-understood geometry---with holonomy strictly contained in $G_{2}$--- and then apply analytical gluing or singular perturbation techniques. Such techniques have been used very widely in global differential geometry over the past three decades and in particular dominate much of the work so far in areas related to exceptional holonomy. There is also a contrast between the two constructions. As the parameter $\epsilon$ tends to $0$ the diameters in the Joyce examples are bounded but the maximum of the curvature tends to infinity; while in the Kovalev examples, as $L\rightarrow \infty$ the curvature stays bounded but not the diameter.

\subsection{The Hitchin functional, moduli and general existence questions.}

The condition for $G_{2}$ holonomy comes in two parts $d\phi=0, d*_{\phi} \phi=0$ and so we may consider the weaker condition of a {\it closed positive 3-form}. This structure is somewhat analogous to a symplectic structure: it is a closed differential form satisfying an open condition pointwise. In this framework, the $G_{2}$-holonomy condition has a variational formulation, due to Hitchin \cite{kn:H1}, \cite{kn:H2}.
We first return to the elementary geometry of 3-forms on a  7-dimensional
vector space $V$. A positive 3-form $\phi$ defines a volume form $\nu(\phi)\in
\Lambda^{7} V^{*}$. In this setting, the 4-form $*_{\phi}\phi$ appears as the
derivative of the function $\nu$:
\begin{equation}  \delta \nu =  \frac{1}{3} *_{\phi}\phi \wedge (\delta \phi). \end{equation} Now let $M$ be a compact oriented $7$-manifold. For   $c\in H^{3}(M;\bR)$ write ${\cal P}_{c}$ for the set of closed positive $3$-forms representing $c$ (of course this could be the empty set). Hitchin considers the functional on ${\cal P}_{c}$ given by the total volume:

$$  V(\phi) =\int_{M} \nu(\phi) . $$
A variation within ${\cal P}_{c}$ is give by the exterior derivative of a 2-form: $\delta \phi=d\alpha$ and after applying (5) and integration by parts we have
$$  \delta V = \frac{1}{3}\int_{M} \alpha \wedge d(*_{\phi}\phi)), $$
so the condition $d*_{\phi}\phi=0$ is the Euler-Lagrange equation defining the critical points of the volume functional $V$.
By  Hodge theory considerations, Hitchin showed that any critical point of is a strict local maximum, modulo the action of the diffeomorphisms of $M$. That is, the Hessian of $V$ is negative-definite transverse to the orbits of the diffeomorphism group acting on ${\cal P}_{c}$.

This approach leads to a simple treatment of the deformation theory of $G_{2}$-structures (which was analysed first by Bryant). The nondegeneracy of the Hessian means that for $c'$ close to $c$ in $H^{3}(M;\bR)$ there is a small perturbation of the critical point in ${\cal P}_{c}$ to one in ${\cal P}_{c'}$, and this is unique up to diffeomorphisms. Let ${\rm Diff}_{0}$ be the identity component of the diffeomorphism group of $M$ and ${\cal T}$ the \lq\lq moduli space'' of $G_{2}$-structures on $M$ modulo  ${\rm Diff}_{0}$. Then the discussion above leads to the fact that the period map $\phi\mapsto [\phi]$ induces a local homeomorphism from ${\cal T}$ to $H^{3}(M;\bR)$.

This approach also gives motivation for the {\it Laplacian flow} on closed positive 3-forms, which was introduced by Bryant \cite{kn:Br3}. This the evolution equation
\begin{equation}   \frac{\partial \phi}{\partial t}= \Delta_{\phi} \phi, \end{equation}
where $\Delta_{\phi}$ is the Laplacian on the metric $g_{\phi}$ defined by $\phi$. This flow has a similar character to the Ricci flow: in fact under the flow (6) the metric evolves by
$$    \frac{\partial g}{\partial t}= -2{\rm Ricci} + Q( \tau), $$
where $Q(\tau)$  is a quadratic expression in the torsion tensor  $\tau=d *_{\phi}\phi$. The leading term (involving two derivatives of $\phi$) is the same as the Ricci flow.  Under the Laplacian flow the volume evolves as
$$   \frac{d V}{dt} = \Vert d *_{\phi}\phi\Vert^{2}, $$
and the flow can be viewed as the ascending gradient flow of the volume functional.

This variational point of view leads to many questions. The most optimistic, naive, hope might be that there is a unique  $G_{2}$-structure (up to diffeomorphism) for each connected component of ${\cal P}_{c}$; that this   is a global maximum of the volume functional and that the Laplacian flow starting from any initial point  exists for all time and converges to this maximum.  If all this were true the existence problem for $G_{2}$-structures would essentially be reduced to understanding the existence of closed positive 3-forms. Examples of Fern\'andez show that this most naive picture does not always hold (\cite{kn:F}, and see also the discussion in \cite{kn:Br3}). Fern\'andez considers left-invariant structures on  a $7$-dimensional nilpotent Lie groups, which descend to compact quotient manifolds $N/\Gamma$. She gives an example where:
\begin{itemize}
\item There are closed positive $3$-forms but no $G_{2}$ structures;
\item The volume functional is not bounded above;
\item The flow exists for all time but the volume tends to infinity and the Riemannian curvature tends to $0$. 
\end{itemize}
But it remains possible that some modified version of the naive picture is true---for example it could (as far as the author knows) be true for manifolds $M$ with $H^{1}(M;\bR)=0$. Or it could be that one needs to allow singularities in the flow and the limits. But at present there is very little known about any of these questions, both for the existence of closed positive 3-forms and for their relation to the existence of $G_{2}$-structures. (For example, it is not known---as far as the author is aware---if there is a closed positive 3-form on the sphere $S^{7}$.)

There is a symmetry in the differential geometric foundations between the 3-form $\phi$ and the 4-form $*_{\phi}\phi$. One can equally well start with a \lq\lq positive''  4-form and use that to define the 3-form. Thus we can ask analogous questions on the relation between \lq\lq co-closed'' 3-forms  and $G_{2}$-structures. By contrast to the closed case, the existence of co-closed structures is completely understood through work of Crowley and Nordstr\"om \cite{kn:CrowNord}. They show that these structures obey an \lq\lq h-principle'', so the question is  reduced to the  homotopy theory of reductions of the structure group of the tangent bundle to $G_{2}$,  which we will discuss in the next section.

\subsection{Some recent developments}

We will now outline some developments in the study of  $G_{2}$-structures from the past 5 years or so. One such development came in the work of Corti, Haskins, Nordstr\"om and  Pacini  \cite{kn:CHNP}, which  clarified and greatly extended Kovalev's work on twisted connected sums. As we sketched in Section 3.2 above, Kovalev's construction begins with \lq\lq building blocks'' which are 3-dimensional Fano manifolds. Recall that a complex manifold is Fano if its anticanonical line bundle $K_{Z}^{-1}$ is ample. Corti {\it et al} considered more generally {\it semi-Fano} 3-folds, where $K_{Z}^{-1}$ is {\it big} and {\it nef}, and showed that these could be used in the twisted connected sum construction. Fano 3-folds, up to deformation, are completely classified and there are 105 deformation types. While there is no complete classification, the number of deformation types of semi-Fano 3-folds is known to be many orders of magnitude larger and so the work of Corti {\it et al} lead to a huge increase in the number of known examples of compact manifolds with  $G_{2}$-structures. The original work of Joyce lead to about 250 examples and Kovalev's original construction gave a  number of a similar order. Corti et al found at least 50 million---and probably many, many more---deformation classes of matching pairs which could be used to construct $G_{2}$-manifolds. They were also able to adjust the construction to give examples with various interesting features and they gave a detailed analysis of the topology of the manifolds constructed.

The topological theory is relatively straightforward in the case of 2-connected 7-manifolds $M$ with $H_{3}(M;\bZ)$ torsion-free. Up to homeomorphism the only invariants are the third Betti number $b_{3}$ and the divisibiity of the Pontrayagin class $p_{1}(M)\in H^{4}(M;\bZ)$. Any such manifold is homeomorphic to a connected sum
$$   \widetilde{S^{3}\times S^{4}}\ \sharp\  S^{3}\times S^{4}\ \dots \sharp \ S^{3}\times S^{4}, $$
where the factor $ \widetilde{S^{3}\times S^{4}}$ (which carries the Pontryagin class) is one of standard list of $S^{3}$ bundles over $S^{4}$. This also gives the diffeomorphism classification up to connected sums with one of the 27 exotic 7-spheres. Corti {\it et al} show that many hundreds of manifolds of this kind can be constructed as twisted connected sums. Moreover they show that the same differentiable manifold can often be constructed in a many different ways, leading to the question whether the resulting $G_{2}$-structures are in the asme connected component of the moduli space.

Topological aspects of   $G_{2}$-manifold theory were developed in a series of papers of Crowley, Nordstr\"om \cite{kn:CrowNord}\cite{kn:CrowNord2} and Crowley, Goette, Nordstr\"om \cite{kn:CGN}. The quotient ${\rm Spin}(7)/G_{2}$ is a 7-sphere and it follows immediately that any spin 7-manifold  has a compatible topological $G_{2}$-structure (in the sense of reduction of the structure group of the tangent bundle). But there are different homotopy classes of such reductions. For example if the tangent bundle is trivial the reduction corresponds to a map
$$  M\rightarrow {\rm Spin}(7)/G_{2}=S^{7}$$
which has an integer degree. Taking account of the action of the diffeomorphisms, Crowley and Nordstrom showed that this integer is cut down to an invariant $\nu\in \bZ/48$. If the divisibility of the Pontrayagin class divides 224 they show that this is a complete invariant (up to homotopy and diffeomorphism), and in general they introduced another more refined invariant to complete the classification. The $\nu$ invariant can be defined as the reduction mod 48 of
$     \chi(W)-3\sigma(W)$,  where $W$ is an 8-manifold with boundary $M$ and with a ${\rm Spin}(7)$ structure compatible with the $G_{2}$ structure on the boundary. They show that $\nu=24$ for all the structures obtained from twisted connected sums, but that there is an example from the Joyce construction with $\nu$ odd.  In \cite{kn:CGN}, Crowley, Goette and Nordstr\"om develop an analytical approach based on the Atiyah, Patodi, Singer theory and define a $\bZ$-valued  lift  $\hat{\nu}$ of $\nu$ which is a deformation invariant of  $G_{2}$-structures. The arguments hinge on the fact that the kernel of the Dirac operator defined by such a structure has fixed dimension $1$, spanned by the covariant constant spinor field. They apply this to the question raised at the end of the previous paragraph, giving an example of a pair of $G_{2}$-structures on the same differentiable manifold which are equivalent at the homotopy level but which can not be joined by a path of torsion-free structures.  In another direction, Crowley and Nordstr\"om define a \lq\lq generalised Eels-Kuiper invariant'' of the smooth structure,  and show that different twisted connected sums can realise different values of this invariant on the same topological manifold \cite{kn:CrowNord2}.

Another line of activity makes contact with Riemannian convergence theory.
 As we noted in 3.2 above, the Joyce and Kovalev constructions can both be seen as describing regions near the \lq\lq boundary of moduli space''---or in other words different kinds of degeneration of sequences of torsion-free $G_{2}$-structures. Understanding all such degenerations is an important general question.   Since these structures define metrics with zero Ricci curvature the substantial body of Riemannian convergence theory, as developed by Cheeger, Colding and others,  can be brought to bear on this question and fit them into a wider context. We recall the basic Gromov compactness theorem for a sequence of complete Riemannian m-manifolds $(M_{i},g_{i})$ with (say) zero   Ricci curvature and with chosen  base points $p_{i}\in M_{i}$. After passing to a subsequence there is a Gromov-Hausdorff limit  $(M_{\infty}, p_{\infty})$ which is a based metric space.  If the sequence has bounded diameter then there is no need to choose base points. In this theory there is a fundamental distinction between {\it non-collapsing} and {\it collapsing} cases.  In the non-collapsing  case the limit is a Riemannian $m$-manifold outside a singular set of codimension at least 4. This occurs for Joyce's examples, when the parameter $\epsilon$ tends to zero and the limit is an orbifold. The twisted sum construction, with parameter $L\rightarrow \infty$, is also non-collapsing, but there are three possible choices of limits depending on the choice of base-points.  In a collapsing situation the Gromov-Hausdorff limit is a space of dimension less than $m$ and one can hope that it can be endowed with some vestige of the structure on the sequence $M_{i}$.

In \cite{kn:D3} the author proposed a programme to study collapsing of  $G_{2}$-structures along co-associative fibrations. This develops ideas of Kovalev \cite{kn:Kov2} and Baraglia \cite{kn:Bar},  and is in a similar vein to seminal work of  Gross and Wilson in the case of $K3$ surfaces collapsing along elliptic fibrations \cite{kn:GW}. Recall that a 4-dimensional submanifold $X$ of a $G_{2}$ manifold $(M,\phi)$ is called co-associative if the restriction of $\phi$ to $X$ vanishes. We consider a smooth map $\pi:M\rightarrow B$ where $B$ is a 3-manifold and $\pi$ is a fibration outside a link $L\subset B$, with fibres diffeomorphic to the K3 surface.  Transverse to $L$, the map is required to be modelled on a non-degenerate complex quadratic form. so the fibres over $L$ are modelled on complex surfaces with ordinary double points. A particular example of this set-up is when $M= Y\times S^{1}, B=S^{2}\times S^{1}$ for a Calabi-Yau 3-fold $Y$ with a Lefschetz fibration over $S^{2}$. The proposal of \cite{kn:D3} is  a sequence of structures of this kind in which the fibres  collapse should have  an \lq\lq adiabatic limit''.  Locally, on a small ball $B_{0}\subset B\setminus L$, the limiting data is given by a family of hyperk\"ahler sructures on a fixed $4$-manifold $X$ (diffeomorphic to a K3 surface), parametrised by $B_{0}$. By the Torelli Theorem for K3 surfaces
this data is a map $$  f:B_{0}\rightarrow {\rm Gr}_{3}^{+}(H^{2}(X,\bR)), $$ where the target  is the space of positive 3-dimensional subspaces of $H^{2}(X)$, the latter being endowed with its cup product form of signature $(3,19)$. The proposal is that the vestige of the $G_{2}$-structures is the requirement that $f$ be the Gauss map of a {\it maximal submanifold} in $H^{2}(X;\bR)$ (i.e. a critical point of the induced volume functional on 3-dimensional submanifolds). This local discussion can be extended to the whole of $B$ but the set-up is more complicated and brings in the monodromy action on the cohomology of the fibres. The main analytical difficulty in carrying through this proposal involves the behaviour around the singular fibres. Indeed just the same difficulty arises when studying the analogous question for  fibred Calabi-Yau 3-folds (a problem which has been studied by many authors, for example \cite{kn:GTZ}).   Significant progress in this direction comes from recent work of Y.Li \cite{kn:Yang},
who constructed a new Ricci-flat metric on $\bC^{3}$ which makes a likely model for the collapsing sequence. Li's work was extended soon after by Conlon and Rochon \cite{kn:CR} and Sz\'ekelyhdi \cite{kn:Gabor3}, leading to many new families of non-compact Ricci-flat K\"ahler metrics. One can expect that these metrics will be important both in K\"ahler geometry and in exceptional holonomy.

Other recent work involving collapsing comes in work of Foscolo, Haskins and Nordstrom \cite{kn:FHN} on codimension one collapse. They consider a circle bundle $\pi:M\rightarrow Y$ over a non-compact Calabi-Yau 3-fold $Y$ and $S^{1}$-invariant  $G_{2}$-structures on $M$. There is a parameter $\epsilon$ related to  the length of the circle fibres and they obtain a well-defined adiabatic limit of the equations as $\epsilon\rightarrow 0$. They give general conditions under which this limiting equation can be solved and show that one can go back to construct a solution of the original problem, for small $\epsilon$. They also find explicit examples of Calabi-Yau 3-folds satisfying the  conditions, so the result is that they get new examples of families of $G_{2}$-structures, collapsing to a 6-dimensional limit. The existence of the relevant metrics on $Y$ depends on recent developments in the K\"ahler theory, of the kind discussed in the first part of this article.

Further developments have taken place in the study of the Laplacian flow  for closed positive 3-forms. Lotay and Wei develop may foundational results and show that the flow can be continued for all time provided that the torsion and curvature are bounded \cite{kn:LotWei}. G. Chen showed that the deep results of Perelman for Ricci flow can be adapted to the Laplacian flow \cite{kn:GChen}. This was used by Fine and Yao to obtain long time existence of a 4-dimensional reduction of the flow under the assumption that the scalar curvature remains bounded \cite{kn:FY}.

 \subsection{Gauge Theory and Calibrated geometry}

 There are two particularly interesting classes of differential geometric objects which can be considered on manifolds with the exceptional holonomy groups $G_{2}, {\rm Spin}(7)$:
\begin{itemize} \item special \lq\lq instanton'' solutions of the Yang-Mills equation;
\item special \lq\lq calibrated'' solutions of the minimal submanifold equation.
\end{itemize}
In the Yang-Mills case, for a Riemannian manifold $M$ with holonomy group $G$ we get a sub-bundle $\Lambda^{2}_{G}$ of the $2$-forms corresponding to the Lie algebra of G (that is, $\Lambda^{2}$ can be viewed as the Lie algebra of the orthogonal group of the tangent bundle and $G$ is a subgroup of that orthogonal group). So we can define an \lq\lq instanton'' to be a connection on a bundle $E\rightarrow M$ whose curvature lies in $\Lambda^{2}_{G}\otimes {\rm End} E $. For general holonomy groups $G$ this may not be an interesting notion, but in the two exceptional cases it yields an elliptic equation (suitably interpreted) for the connection. In \cite{kn:DT} the author and R. Thomas pointed out the strong formal analogies between these equations in dimensions 8 and 7 and the usual instanton theory and Floer theory in   dimensions 4 and 3. In fact the discussion can be extended to link dimensions 8,7,6 and the holonomy groups ${\rm Spin}(7), G_{2}, SU(3)$. The \lq\lq instanton'' equation over a manifold with holonomy $SU(3)$ is just the Hermitian Yang-Mills equation discussed in Section 2.2 above, whose solutions (over a compact manifold) correspond to stable holomorphic vector bundles. 

In the submanifold case we have the calibrated submanifolds:
\begin{itemize} \item ${\rm Spin}(7)$: codimension-4 {\it Cayley submanifolds};
\item $G_{2}$: codimension 4 {\it associative submanifolds}, codimension 3 {\it co-associative submanifolds}
\item $SU(3)$: codimension 4 {\it holomorphic curves}, codimension 3 {\it special Lagrangian submanifolds}
\end{itemize}
The calibrated condition means that these submanifolds are absolutely volume minimising in their homology class. This theory goes back to the seminal paper of Harvey and Lawson \cite{kn:HL} from the early 1980's  (featured in Harvey's ICM 1983 contribution \cite{kn:H}) which was influential in reviving interest generally in exceptional holonomy. A notable feature of the cases above is that the calibrated condition is an elliptic equation, suitably interpreted. 
The suggestion in \cite{kn:DT} was that there could be \lq\lq enumerative'' theories built on these geometric objects (connections and submanifolds), yielding numerical invariants or other structures such as Floer homology groups. In the two decades since \cite{kn:DT} only one such theory has been developed rigorously, leading to what are sometimes called Donaldson-Thomas invariants  of Calabi-Yau 3-folds. This development has been entirely algebro-geometric, working with stable holomorphic bundles, or more generally sheaves, rather than Hermitian Yang-Mills connections. In fact the main algebro-geometric interest has been in the rank 1 case,  where the Hermitian Yang-Mills theory does not directly apply.

Going back to differential geometry: there has been substantial activity and advances. A first question is the very existence of these higher dimensional instantons. This has been addressed by Walpuski \cite{kn:Walp1},\cite{kn:Walp2} S\'a Earp and Menet \cite{kn:SaW}, \cite{kn:MNE} through gluing constructions which parallel the Joyce and Kovalev constructions of the underlying manifolds.  S\'a Earp showed that a stable holomorphic bundle over a Fano building block whose restriction to the divisor $D$ is stable gives rise to a Hermitian -Yang-Mills connection on the cylindrical-end manifold.  Two such bundles can be used to construct an instanton the twisted connected sum if they satisfy a suitable matching condition over the divisors. In \cite{kn:Walp1}, Walpuski proved the existence of instantons over one of Joyce's examples by gluing a connnection constructed from a 4-dimensional instanton over the Eguchi-Hanson manifold to a flat connection over the orbifold. Walpsuki used related techniques to construct instantons over manifolds with ${\rm Spin}(7)$-holonomy \cite{kn:Walp2}.

The fundamental problem in developing enumerative theories  is the potential failure of compactness. The foundations to treat this question were developed by Tian \cite{kn:Ti2}, leading to a theory which has close analogies with the Riemannian convergence theory.  For any sequence of solutions of the instanton equation on a fixed bundle there is a subsequence which converges outside a set $S$ of codimension at least $4$ and the  codimension 4 component of $S$ satisfies a generalised form of the corresponding calibrated condition. Restricting to the $G_{2}$ case, the basic problem occurs in a generic 1-parameter family $\phi_{t}$ of $G_{2}$-structures. We expect that for each $t$ there are a finite number of $G_{2}$ instantons (on a  given bundle $E$), but at some parameter value $t_{0}$ this number may change because a family of solutions $A_{t}$ for $t<t_{0}$ (say) may diverge as $t\rightarrow t_{0}$ with \lq\lq bubbling'' along some associative submanifold $P$. Thus, starting on the gauge theory
side, we are inevitably drawn to consider the submanifold theory as well.
This fits in with the close analogies between the theories at the formal
level and also with the algebro-geometric point of view in the Calabi-Yau case (where we asociate an ideal sheaf to a curve). Under reasonable assumptions  one can show that this divergent phenomenon occurs exactly for the parameter values $t$ where there is a solution of a certain \lq\lq Fueter equation'' over an associative submanifold $P$ \cite{kn:Walp3}.

 Haydys and Walpuski have proposed an approach to overcome this fundamental compactness problem \cite{kn:HayWalp}. They consider a variant of the Seiberg-Witten equations over a 3-manifold and their proposal is to count both instantons  and associative submanifolds, the latter weighted by a count of the solutions of the Seiberg-Witten equations. The reason for this is that compactness can also fail for these Seiberg-Witten equations and an analysis of this failure of compactness suggests that it occurs exactly when there is a solution of the same Fueter equation. So one can hope that (with a careful choice of signs) the jump in the count of instantons should be exactly compensated by the jump in the Seiberg-Witten count.

There are many technical problems in carrying through this Haydys and Walpuski programme (particulary in the case when $P$ is a homology sphere and reducible solutions present difficulties, of a kind familiar in the gauge theory literature). The programme motivates a better understanding of higher codimension singularities, making contact with deep work of Taubes \cite{kn:Taubes}. As shown by Doan and Walpuski \cite{kn:DW}, one has to consider also solutions of the Fueter equation which have singularities along 1-dimensional sets in $P$ and it seems likely that these correspond to instantons on the 7-manifold with $1$-dimensional singular sets. Once again there is link here with the complex geometry. Bando and Siu showed that the Hermitian-Yang-Mills theory can be extended to reflexive sheaves and one expects that, transverse to the 1-dimensional singular set in the 7-manifold, the geometry of the connection should be modelled by a reflexive sheaf on a complex 3-fold, with a point singularity.


\end{document}